\documentclass[a4paper,12pt]{amsart}

\usepackage{amsthm}
\usepackage{amsmath}
\usepackage{amsfonts}
\usepackage{amssymb}
\usepackage{amsopn}
\usepackage{amscd}
\usepackage{graphicx}
\usepackage{enumerate}
\usepackage[all,cmtip]{xy}
\usepackage{url}
\usepackage{color}
\theoremstyle{plain}

\newtheorem*{thm-main}{Main Theorem}

\theoremstyle{definition}

\theoremstyle{remark}

\DeclareMathOperator{\N}{\mathbf{N}}

\begin{document}

\author[J.~Byszewski]{Jakub Byszewski}
\address[J.~Byszewski]{Institute of Mathematics, Jagiellonian University,
{\L}ojasiewicza~6, 30-348 Krak\'ow, Poland}
\email{jakub.byszewski@uj.edu.pl}

\author[F.~Falniowski]{Fryderyk Falniowski}
\address[F.~Falniowski]{Department of Mathematics, Cracow University of Economics,
Rakowicka~27, 31-510 Krak\'ow, Poland}
\email{fryderyk.falniowski@uek.krakow.pl}

\author[D.~Kwietniak]{Dominik Kwietniak}
\address[D.~Kwietniak]{Institute of Mathematics, Jagiellonian University,
{\L}ojasiewicza~6, 30-348 Krak\'ow, Poland}
\email{dominik.kwietniak@uj.edu.pl}

%\title{The Hoehn-Mouron map which is a transitive dendrite map with zero entropy}
\title{Transitive dendrite map with zero entropy}

\begin{abstract}
%L.~
Hoehn and %Ch.~
Mouron [Ergod. Th. \& Dynam. Sys. (2014) \textbf{34}, 1897--1913] constructed a map on the universal dendrite that is topologically weakly mixing but not mixing. We modify the Hoehn-Mouron example to show
that there exists a transitive (even weakly mixing) dendrite map with zero topological entropy. This answers the question of Baldwin [Topology (2001) \textbf{40}, 551--569].
\end{abstract}

\maketitle

\section{Introduction}
Blokh \cite{Blokh} proved that a transitive map of a connected topological graph is either non-invertible and has the relative specification property, or is an irrational rotation of the circle. In the former case, a transitive map has dense set of periodic points and positive topological entropy. There are other connections between transitivity and topological entropy for one-dimensional maps (see \cite{AKLS}).
This suggests a question: what topological conditions force an analogous result? A natural candidate for the family of continua with similar properties is the class of dendrites. Dendrites generalize topological trees and appear in dynamics (for example as Julia sets for some complex polynomials \cite[Thm. 4.2]{CG}). It is also known that maps on dendrites share some dynamical properties with graph maps (see, for example, \cite{AHNO, K, KKM, Mai, N, S}).
In particular, if a dendrite contains a free arc then a transitive map necessarily has positive topological entropy (see \cite{DSS}).
%To some extent it is true (if dendrite contains a free arc) but in general...

%On the other side dendroids are the generalization of dendrites and ... Cantor fan.

%This answers a question of
%Baldwin from 2001

In \cite{HM}, Hoehn and Mouron introduced a family of self-maps on the Wa\.{z}ewski universal dendrite $D_\omega$ and showed that among these maps there are  topologically weakly mixing but not mixing ones. Somewhat imprecisely, we say that a map is a \emph{Hoehn-Mouron map} if it is defined following a similar procedure as in \cite{HM}. %Our main result proves existence of a Hoehn-Mouron map.
%Here we give a self-contained exposition of some variant of the Hoehn-Mouron construction.
%This leads to a proof
This note contains a proof of the following theorem (here we present only the main part of it, for the precise statement see p. \pageref{Main}):

%In what follows we use the notation and facts introduced in \cite{HM}, so the reader should have a copy of that paper at hand.
%\emph{there exists a transitive dendrite map with entropy zero}.

\begin{thm-main}[short version]  %Let $D_\omega$ denote the Wa\.{z}ewski universal dendrite.
There exists a weakly mixing, not mixing, proximal Hoehn-Mouron map $f\colon D_\omega\to D_\omega$ with zero topological entropy.
\end{thm-main}

Our result gives an affirmative answer to the question of Baldwin \cite{Baldwin}, whether there is
a transitive dendrite map with zero entropy. This also provides a new example of a transitive map of a dendrite which
does not have dense periodic points (because proximality excludes existence of more than one minimal point and every periodic point is minimal). This property of the Hoehn-Mouron example is also noted in \cite{AHNO}.
Observe that dendrites have the fixed point
property, hence every continuous dendrite map has at least one fixed point. An example of a transitive dendrite map with
only one fixed point and no other periodic points is given in \cite{S}, but the map constructed there has positive topological entropy. %Another recent work in this area is \cite{AHNO}.

%\section{Main theorem}

%so this note contains only the proof of the main result:

\section{Basic definitions}
We assume the reader is familiar with elementary theory of dynamical systems.

A \emph{dendrite} is a compact, pathwise connected, locally connected, non-separating, and nowhere dense subset of the plane, or equivalently a locally connected continuum that contains no simple closed curves. The \emph{universal dendrite} (the \emph{Wa\.zewski dendrite}, see \cite[\S 10.37, p. 181-185]{Nadler}) is a dendrite $D_\omega$ such that each ramification point of $D_\omega$ is of infinite order and for each arc $A\subset D_\omega$, the set of ramification points of $D_\omega$ contained in $A$ is dense in $A$. Let $f\colon D_\omega\to D_\omega$ be a continuous map.
%Fix a metric $\rho$ compatible with the topology of $D_\omega$.
Given non-empty open sets $U,V\subset D_\omega$, define  $N(U,V)=\{n\in\N:f^{-n}(U)\cap V\neq\emptyset\}$.
We recall that $f$ is \emph{transitive} if $N(U,V)\neq\emptyset$ for every non-empty open sets $U,V\subset D_\omega$.
A map $f$ is \emph{weakly mixing} if the product map $f\times f$ is transitive.
We say that $f$ is mixing if the set $N(U,V)\neq\emptyset$ has finite complement in $\N$ for every non-empty open sets $U,V\subset D_\omega$.
Equivalently, $f$ is weakly mixing if $N(U,V)$ is thick for all nonempty open sets $U,V\subset X$. Recall that a set $A\subset \N$ is
\emph{thick}, if it contains arbitrarily long blocks of consecutive integers, that is,
for every $n>0$ there is $k\in\N$ such that $\{k,k+1,\ldots,k+n-1\}\subset A$. An $f$-invariant closed nonempty set $M\subset D_\omega$ is \emph{minimal} if an orbit of every point $x\in M$ is dense in $M$. A point $x\in D_\omega$ is a \emph{minimal point} if the closure of its orbit is a minimal set. We say that $f$ is \emph{proximal} if it has a unique minimal point. Note that if $x\in D_\omega$ is a unique minimal point for $f$, then $x$ is fixed by $f$.

\section{The main theorem}

\begin{thm-main} \label{Main} Let $D_\omega$ denote the Wa\.{z}ewski universal dendrite. There exists a set $Z\subset\N$ such that the Hoehn-Mouron map $f_Z\colon D_\omega\to D_\omega$ has the following properties:
\begin{enumerate}
  \item $f_Z$ is weakly mixing, but not mixing;
  \item all transitive points of $f_Z$ are end points of $D_\omega$;
  \item $f_Z$ has a unique fixed point $o$, which is, in fact, a unique minimal point of $f_Z$ (it follows that $f_Z$ is proximal);
  \item $f_Z$ is uniquely ergodic, with the only $f_Z$-invariant Borel probability measure being the Dirac measure concentrated on $o$;
  \item $f_Z$ is Li-Yorke chaotic, but not DC$_2$-chaotic;
  \item $f_Z$ has zero topological entropy.
\end{enumerate}
\end{thm-main}
\begin{figure}
\includegraphics{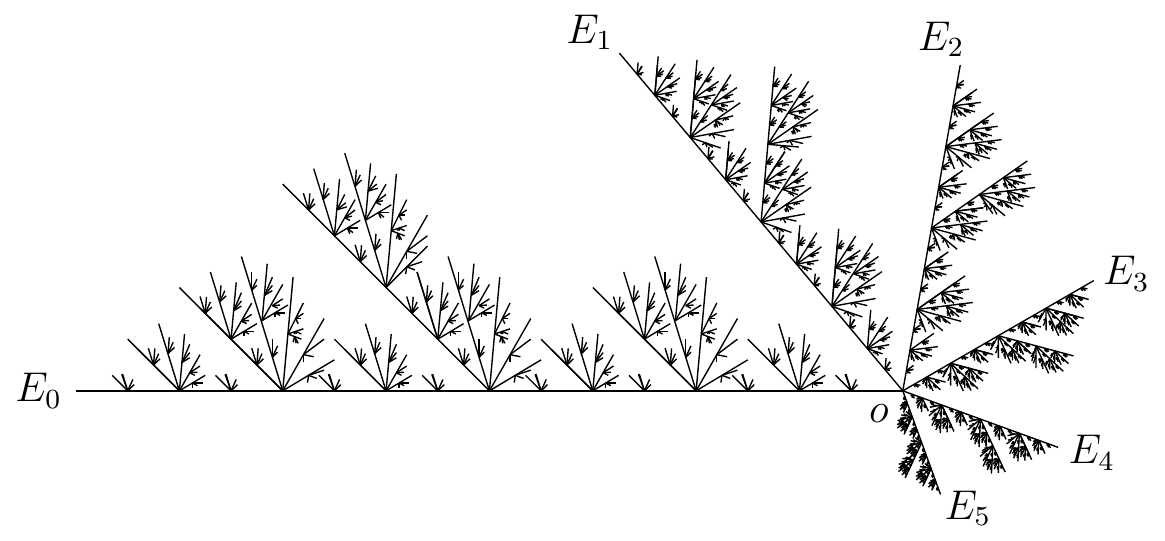}
\caption{An approximation of the Wa\.{z}ewski universal dendrite and subdendrites $E_0,\ldots, E_5$.}
\label{rys}
\end{figure}
\begin{proof}
We follow the notation and terminology introduced in \cite{HM}. We write
\[
D_\omega=\bigcup_{j=0}^\infty E_j, \text{ where }E_i \cap E_j=\{o\} \text{ for }i\neq j,
\]
and each $E_j$ is a homeomorphic copy of $D_\omega$. The subdendrites $E_j$ are defined on page 7 of \cite{HM}.
Let $U_j=E_j\setminus\{o\}$, thus $U_j$ is an open subset of $D_\omega$ for $j=0,1,\ldots$.
Example 17 of \cite{HM} proves that given a set $Z\subset \N$
fulfilling certain assumptions there is a weakly mixing, but not mixing map $f_Z\colon D_\omega\to D_\omega$ satisfying
\begin{align}
f_Z(E_j)&=E_{j-1}, \quad\text{for }j\ge 1,\label{claim1}\\
f_Z^{n}(E_0)&\subset \bigcup_{n+j-1\in Z}E_j, \quad\text{for }n=1,2,\ldots.\label{claims1734}
\end{align}
It remains to find a set $Z$ such that the map $f_Z$ exists and has all the properties listed in the main theorem.
We construct the set $Z$ inductively.
Let $N_0=0$ and $M_0=1$. Assume that we have defined $N_K$ and $M_K>N_K$ for some $K\in\N$. We set
$N_{K+1}=(2^{M_K}+1)\cdot M_K$ and we take $M_{K+1}$ large enough so that
\textsc{Claims} 17.2 and 17.5 of \cite{HM}
are satisfied.\footnote{In the notation of \cite{HM}, it means that if $Z\subset\N$ is a thick set such that $[N_{K+1},M_{K+1}]\cap\N\subset Z$, then
for each $0\le i \le K+1$ there are some $m,n\in\N$ such that $\alpha(s)$ is defined for each
$s\in T_{s_i,\,K+1,\,m}$ and $\alpha^n(T_{s_i,\,K+1,\,m})=T_{\emptyset,\,K+1,\,0}$.}

Take
\[
Z=\bigcup_{K=1}^\infty [N_K,M_K]\cap\N.
\]
Then $f_Z$ as defined in \cite{HM} is weakly mixing but not mixing.
%The set $E_0$ will have forward images which are arbitrarily small. Since each subdendrite $E_j$ is mapped onto $E_0$ in $j$ iterates, it follows that any non-empty open set $U$ which is contained in one of the sets $E_j$ will have, by Proposition \ref{claims1734} forward images which are arbitrarily small. Thus, $f_Z$ is not mixing. Weakly mixing comes from...

% (the proof is the same as \cite[page 11]{HM} --- it uses only two assumptions that $M_K$ is chosen for $N_K$ so that Proposition \ref{claims1734} and Theorem \ref{claim175} are satisfied and that $\N\setminus Z$ is thick).

By \eqref{claims1734} for each $n\in\N$ we have
\[
f_Z^{n}(E_0)\subset \bigcup_{n+j-1\in Z}E_j, \quad\text{for }n=1,2,\ldots.
\]
It follows that
\[
f_Z^{n}(U_0)\subset \bigcup_{n+j-1\in Z}U_j\cup\{o\}, \quad\text{for }n=1,2,\ldots.
\]
But with our $Z$ we know that if $N=M_K+1$ for some $K\in\N$, then
\[
f_Z^{N}(U_0)\subset \bigcup_{j=2^{M_K}\cdot M_K}^\infty U_j\cup\{o\},
\]
hence
\[
f_Z^{N+l}(U_0)\subset \bigcup_{j=2^{M_K}\cdot M_K-l}^\infty U_j\cup\{o\}, \quad\text{for }l=1,2,\ldots,2^{M_K}\cdot M_K-1.
\]
Therefore for each $K\in\N$ we have
\begin{equation}\label{eq:disjont}
f_Z^j(U_0) \cap U_0=\emptyset\quad\text{for }j=M_K+1,\ldots,\left(2^{M_K}+1\right)\cdot M_K.
\end{equation}
Let $x\in U_0$ and $N(x,U_0)=\{n\in\N : f_Z^n(x)\in U_0\}$.
It follows from \eqref{eq:disjont} that
\begin{equation}\label{eq:lower_density}
\liminf_{n\to\infty}\frac{|N(x,U_0)\cap\{0,1,\ldots,n-1\}|}{n}=0.
\end{equation}
If $\bar\mu(U_0)>0$ for some $f_Z$-invariant Borel probability measure $\bar\mu$, then without loss of generality we can assume that $\bar\mu$ is ergodic and then \eqref{eq:lower_density} would contradict the pointwise ergodic theorem.
It follows that $\mu(U_0)=0$ for every $f_Z$-invariant Borel probability measure $\mu$. Moreover,  we also have $\mu(U_j)=0$ for every $j\in\N$ since
\[
U_j\subset f_Z^{-j}(U_0), \quad\text{for }j=1,2,\ldots.
\]
(this follows easily from \eqref{claim1}).
Hence the Dirac measure concentrated on $o$ is the only $f_Z$ invariant probability measure. By the Variational Principle \cite{Walters}, $f_Z$ has topological entropy zero.
The statement concerning Li-Yorke chaos follows from $f_Z$ being weakly mixing (see \cite{Iwanik}). The statement about distributional chaos is a consequence of the fact that for each $x,y\in D_\omega$ the pair $(x,y)$ is proximal with upper Banach density one (because the Dirac measure concentrated on $o$ is the only $f_Z$-invariant probability measure, see \cite{LT}). We refer the reader to \cite{Iwanik, LT} for the details and definitions of Li-Yorke and DC$_2$-chaos.
\end{proof}

\section*{Acknowledgements}
We are greatly indebted to Vladim{\'{\i}}r {\v{S}}pitalsk{\'y} for his remarks and fruitful discussion. We also gratefully thank Logan Hoehn and Piotr Oprocha for several helpful comments concerning the subject of the paper. The research of Jakub Byszewski and Dominik Kwietniak was supported by the  National Science Centre (NCN) under grant no. DEC-2012/07/E/ST1/00185.

\end{document}